\documentclass[10pt,reqno]{amsart}
\RequirePackage[OT1]{fontenc}
\RequirePackage{amsthm,amsmath}
\RequirePackage[numbers]{natbib}
\RequirePackage[colorlinks,linkcolor=blue,citecolor=blue,urlcolor=blue]{hyperref}
\usepackage{amssymb}
\usepackage{anysize}
\usepackage{hyperref}
\usepackage{extarrows}
\usepackage{multirow}
\usepackage{natbib}
\usepackage{stmaryrd}
\usepackage{color}
\usepackage{amsmath}
\usepackage{enumitem}
\usepackage{mathtools} 
\usepackage{dsfont} 
\usepackage{comment}
\usepackage{soul}

\numberwithin{equation}{section}
\newtheorem{theorem}{Theorem}[section]
\newtheorem{lemma}[theorem]{Lemma}

\newtheorem{definition}[theorem]{Definition}
\newtheorem{remark}[theorem]{Remark}

\renewcommand{\Re}{\mathrm{Re}\,}
\renewcommand{\Im}{\mathrm{Im}\,}
\newcommand{\E}{{\mathbb E }}

\newcommand{\R}{{\mathbb R }}
\newcommand{\N}{{\mathbb N}}

\renewcommand{\P}{{\mathbb P}}
\newcommand{\C}{{\mathbb C}}
\newcommand{\ii}{\mathrm{i}}
\newcommand{\deq}{\mathrel{\mathop:}=}
\newcommand{\dd}{\mathrm{d}}
\newcommand{\ie}{\emph{i.e., }}
\newcommand{\eg}{\emph{e.g., }}
\newcommand{\cf}{\emph{c.f., }}

\newcommand{\wt}{\widetilde}

\newcommand{\bs}{\boldsymbol}
\renewcommand{\mathbf}[1]{\bs{#1}}

\def\Tr{\mathrm{Tr}}

\def\one{\mathds{1}}

\def\<{\langle}
\def\>{\rangle}

\def\X{\mathcal{X}}
\def\Dim{\Delta {\mathrm{Im}}\,}

\marginsize{25mm}{25mm}{25mm}{26mm}

\begin{document}

	\begin{minipage}{0.85\textwidth}
		\vspace{2.5cm}
	\end{minipage}
	\begin{center}
		\large\bf Small deviation estimates for the largest eigenvalue of Wigner matrices
		
	\end{center}

	\renewcommand{\thefootnote}{\fnsymbol{footnote}}	
	\vspace{0.5cm}
	
	\begin{center}
		\begin{minipage}{1.4\textwidth}

			\begin{minipage}{0.33\textwidth}
				\begin{center}
					L\'aszl\'o Erd\H{o}s\footnotemark[1]\\
					\footnotesize 
					{IST Austria}\\
					{\it lerdos@ist.ac.at}
				\end{center}
			\end{minipage}
			\begin{minipage}{0.33\textwidth}
				\begin{center}
					Yuanyuan Xu\footnotemark[2]\\
					\footnotesize 
					{IST Austria}\\
					{\it yuanyuan.xu@ist.ac.at}
				\end{center}
			\end{minipage}
		\end{minipage}
	\end{center}
	
	\bigskip

	\footnotetext[1]{Partially supported by ERC Advanced Grant "RMTBeyond" No.~101020331.}
	\footnotetext[2]{Supported by ERC Advanced Grant "RMTBeyond" No.~101020331.}

	\renewcommand{\thefootnote}{\fnsymbol{footnote}}	
	
	\vspace{5mm}
	
\begin{center}
	\begin{minipage}{0.83\textwidth}\footnotesize{
			{\bf Abstract.}} We establish  precise right-tail small deviation estimates for the largest eigenvalue of
		real symmetric and complex Hermitian matrices whose entries are independent random variables
		with uniformly bounded moments. The proof relies on a Green function comparison 
		along a continuous interpolating matrix flow for a long time. Less precise estimates are also obtained in the left tail.

	\end{minipage}
\end{center}

	\vspace{5mm}
	
	{\small
		\footnotesize{\noindent\textit{Date}: March 1, 2022}\\
	}
	
	\vspace{2mm}

	\thispagestyle{headings}

\section{Introduction and main result}

The celebrated Tracy-Widom law~\cite{TW94,TW96} describes the asymptotic distribution
of the largest eigenvalue $\lambda_N$ of Gaussian random  matrices as the dimension $N$ of the 
matrix tends to infinity. By  well-established {\it universality} results \cite{AEKS20,BEY14,EYY12,FS10, HLY20,LY17,LS15,LS18,LY14, PY14, S99, TV10} the same limiting distribution holds for a large class of random matrices.
The Tracy-Widom distribution is asymmetric and it has characteristic left and right tails.
To  what extent  do the large deviation behaviour of the largest eigenvalue $\lambda_N$ 
reflect the tail behaviour of the Tracy-Widom tails?

Unsurprisingly, without additional strong decay condition on the distribution of the matrix elements,
the large  $N$ limit is not interchangeable with the large tail asymptotics. But even with additional conditions,
standard techniques developed for the  universality proofs are inadequate to handle the 
large deviation regime and the strongest existing {\it large deviation} results~\cite{A16,AGH21,GH20,GM20,H20,M21}
use very different approaches.
In this work we study the {\it small deviation} regime with the {\it Green function comparison} method,
originally developed for universality proofs. The advantage of our approach is twofold: (i) only  the usual condition on 
finiteness of all  moments for the matrix elements is required, and (ii) we prove a very accurate
tail behavior; in fact, we can transfer  the existing best results for Gaussian matrices  
to general Wigner matrices without any loss of precision. Under such weaker moment conditions,
however, only the regime of  not very  large  deviation can be handled;
in this regime the scale of probability is only inverse polynomial in $N$
instead of the customary (sub)exponentially smallness in standard large deviation theory.
In fact, an even larger regime was coined "small deviation" in \cite{A05} and \cite{L07}.
However, in statistical applications the standard large deviation theory for the largest eigenvalue
of a sample covariance matrix captures the regime where the {\it miss probability} is
exponentially small in the number of samples, see e.g. Theorem 2 in \cite{BDMN111} and Theorem 5.1 in \cite{FVK08},  which may be impractically small.
Our  result enables formulating similar relations for polynomially small miss probabilities which may be more realistic in certain
applications.
On the methodological  side, quite importantly, our result demonstrates that Green function comparison is able to 
handle events of arbitrary polynomially small probability.

\bigskip

Now we describe our setup more precisely, 
first by collecting the Gaussian results and some facts about the Tracy-Widom law.
Consider a real symmetric random matrix $H$ of dimension $N$, 
whose entries are independent Gaussian random variables up to symmetry, \ie for any $1\leq i\leq j \leq N$,
\begin{align}\label{GOE}
	\sqrt{N}H_{ij} \stackrel{{\rm d}}{=} \mathcal{N}(0,1), \qquad \sqrt{N}H_{ii} \stackrel{{\rm d}}{=} \mathcal{N}(0,2), \qquad H_{ij}=H_{ji}.
\end{align}

The resulting random matrix model is called the Gaussian Orthogonal Ensemble (GOE). The complex Hermitian version is called the Gaussian Unitary Ensemble (GUE), where the matrix entries are complex-valued Gaussian random variables with independent real and imaginary parts, \ie for any $1\leq i\leq j \leq N$,
\begin{align}\label{GUE}
	\sqrt{N} \Re H_{ij} \stackrel{{\rm d}}{=} \mathcal{N}(0,1/2), \quad \sqrt{N} \Im H_{ij} \stackrel{{\rm d}}{=} \mathcal{N}(0,1/2),  \quad H_{ij}=\overline{H_{ji}}, \quad \sqrt{N}H_{ii} \stackrel{{\rm d}}{=} \mathcal{N}(0,1).
\end{align}
We will use the parameter $\beta$ to indicate the symmetry class, \ie $\beta=1$ for real symmetric 
and $\beta=2$ for complex Hermitian matrices. We also refer to these Gaussian ensembles as G$\beta$E for short.

The eigenvalues of the matrix $H$ are denoted by $(\lambda_j)_{j=1}^N$ labelled in a non-decreasing order. 
The limiting distributions of the normalized fluctuations of the largest eigenvalue $\lambda_N$ were identified 
by Tracy and Widom in \cite{TW94,TW96}, \ie for any fixed $x\in \R$,  
\begin{align}\label{le one}
	\lim_{N \rightarrow \infty} \P^{\mathrm{G \beta E}} \Big( N^{2/3} (\lambda_N-2)\leq x \Big) =\mathrm{TW}_{\beta}(x)\,.
\end{align}
The limiting distribution functions 
$\mathrm{TW}_{\beta}$
are given by
\begin{align}
	\mathrm{TW}_2(x):=\mathrm{e}^{-\int_x^{\infty} (t-x) q^2(t) \dd t}, \qquad \mathrm{TW}_1(x):=\sqrt{\mathrm{TW}_2(x)} \mathrm{e}^{-\frac{1}{2} \int_x^\infty q(t) \dd t},
\end{align}
where 
$q=q(x)$  is the solution to the Painlev\'e II equation with asymptotics given by the Airy function:
\begin{align}
	q''(x)=xq(x)+2 q^3(x), \qquad q(x) \sim \mathrm{Ai}(x) \sim \frac{e^{-\frac{2}{3} x^{3/2}}}{2 \sqrt{\pi} x^{1/4}}, \quad x \rightarrow \infty. 
\end{align} 
The right and left tail asymptotics of $\mathrm{TW}_\beta(x)$ are well known (see, e.g. Chapter 3 \cite{AGZ10})
\begin{equation}\label{TW_right}
	1-\mathrm{TW}_\beta(x) \sim x^{-\frac{3 \beta}{4}}  e^{-\frac{2\beta}{3} x ^{3/2}}, \quad x \rightarrow +\infty; \qquad \mathrm{TW}_\beta(x) \sim |x|^{-\frac{\beta}{16}}e^{-\frac{\beta}{24} |x|^{3}}, \quad x \rightarrow -\infty.
\end{equation}

Inspired by (\ref{le one}), similar asymptotics of the Tracy-Widom distributions in (\ref{TW_right}) are expected to be valid for the tail distribution of the largest eigenvalue, at least if $x$ is moderately large. 
In this paper, we focus on the {\it small deviation regime}, \ie $x=O((\log N)^{2/3})$ for the right tail and $-x=O((\log N)^{1/3})$ for the left tail and study the precise asymptotics as indicated in (\ref{TW_right}). 
Much larger $x$ regimes, but still well below the standard large deviation regime with exponential decay in $N$, were also considered in \cite{PZ17} to prove a law of iterated logarithm for the largest eigenvalue of the GUE. 
Here a precise right-tail estimate for the largest eigenvalue was obtained, \ie there are constants $C>0$ and $\delta>0$ so that for any $1\leq x\leq \delta N^{1/6}$, 
\begin{align}\label{gue_sharp}
	C^{-1} x^{-3/2} e^{-\frac{4}{3} x^{3/2}} \leq \P^{\mathrm{GUE}} \Big( N^{2/3} (\lambda_N-2) \geq x \Big) \leq Cx^{-3/2} e^{-\frac{4}{3} x^{3/2}}.
\end{align}
Before that, less precise estimates for the tail distributions of the largest eigenvalue of the Gaussian ensembles were obtained in \cite{LR10} for an even  larger range of $x$, see also \cite{A05} for the 
GUE and \cite{L07} for the GOE. In particular, for any $0<x \leq N^{2/3}$, there exists a numerical constant $C_0>0$ such that
\begin{align}\label{less_precise}
	C_0^{-1} e ^{-C_0\beta x^{3/2}} \leq 	\P^{\mathrm{G\beta E}} \Big( N^{2/3} (\lambda_N-2) \geq x \Big) \leq C_0 e ^{-\beta x^{3/2}/C_0}; \nonumber\\ 	 
	C_0^{-1} e ^{-C_0\beta x^{3}} \leq \P^{\mathrm{G\beta E}} \Big( N^{2/3} (\lambda_N-2) \leq -x \Big) \leq C_0 e ^{-\beta x^{3}/C_0}.
\end{align}
We note that the sharp right tail asymptotics~\eqref{gue_sharp} has apparently been obtained only in the GUE case
and not for the GOE. For completeness,  using the uniform Plancherel-Rotach asymptotic estimates for Hermite polynomials \cite{S59},
in the Appendix we give the proof of the following  
precise upper bound for the right tail of the GOE: 
\begin{lemma}\label{lemma_tail}
	For any $1\leq x\leq N^{1/4}$, there exists some constant $C>0$ such that
	\begin{align}\label{tail}
		\P^{\mathrm{G \beta E}} \Big( N^{2/3} (\lambda_N-2)>x \Big) \leq Cx^{-\frac{3 \beta}{4}} e^{-\frac{2\beta}{3} x ^{3/2}}.
	\end{align}
\end{lemma}
We remark that  estimates for the left-tail distribution of the Gaussian ensembles in general are less accessible and 
may be obtained using refined Riemann-Hilbert analysis as discussed in \cite{PZ17} (\cf \cite{BDMMZ01}). In particular, even the precise constants in the exponents in the second line of~\eqref{less_precise} are not known,
in contrast to~\eqref{gue_sharp}, where  the polynomial prefactor is also precisely identified.

Beyond the Gaussian ensembles, the right-tail upper bound in (\ref{less_precise}) was extended to non-Gaussian matrices \cite{FS10,L09} 
when the matrix entries have sub-Gaussian tails. An exponential upper bound $e^{-cx}$ for $x=O((\log N)^{\log N})$ was obtained in 
Theorem~2.2 of \cite{EYY12} for Wigner matrices whose entries have sub-exponential tails. A standard large deviation principle 
for $x \sim N^{2/3}$ with speed $N$ and an explicit rate function was established in \cite{AGH21,GH20} for Wigner matrices whose entries have sub-Gaussian tails.

In this paper, we show that the right-tail small deviation estimates with precise power as in (\ref{tail}) can be extended to non-Gaussian matrices in the regime up to 
$x=O\big((\log N)^{2/3}\big)$, under a weaker condition that the matrix entries have uniformly bounded moments. 
Less precise left-tail estimates can also be obtained in combination with (\ref{less_precise}) up to  $-x=O\big((\log N)^{1/3}\big)$.

\subsection{Main result}
Let~$H_N$ be a real symmetric or complex Hermitian random matrix of size $N$, whose entries are  independent but not necessarily identically distributed random variables up to symmetry $H_{ij}=\overline{H_{ji}}$. For real symmetric matrices, we assume that for any $1\leq i < j \leq N$, 
\begin{align}\label{second_moment}
	\E[H_{ij}]=0. \qquad \E[ |\sqrt{N} H_{ij}|^2]=1, \qquad \E[H_{ii}]=0,\qquad  \E [|\sqrt{N} H_{ii}|^2]=:m_2 <\infty.
\end{align}
For complex Hermitian matrices, we further assume that 
\begin{align}\label{second_moment_complex}
	\E[(H_{ij})^2]=0, \qquad 1\leq i < j \leq N,
\end{align}
that holds, e.g., if the real and imaginary parts of $H_{ij}$ have the same distribution.
We also assume that all moments of the entries of $\sqrt{N}H_N$ are uniformly bounded, \ie for any $k \geq 3$, there exists $m_k$ independent of $N$ such that
\begin{equation}\label{moment_condition}
	\max_{i,j}\E [|\sqrt{N} H_{ij}|^k] \leq m_k\,.
\end{equation}	
Note that we do not assume identical distribution for $H_{ij}$, so our model is slightly more general
than the customary Wigner ensemble, but for simplicity we will refer to as the Wigner ensemble. 
Our main theorem shows that the largest eigenvalue $\lambda_N$ of $H_N$ has the following small deviation estimates that exactly matches the best known Gaussian results: 
\begin{theorem}\label{main}
	Fix $M>0$. Let $\lambda_N$ be the largest eigenvalue of the real symmetric
	or complex Hermitian $H_N$ satisfying~\eqref{second_moment}--\eqref{moment_condition}. 
	Then for any $1 \leq x \leq M(\log N)^{2/3}$, there exists a constant $C=C(M)>0$ such that
	\begin{equation}\label{tw_tail}
		\P \Big( N^{2/3} (\lambda_N-2)>x \Big) \leq C x^{-\frac{3 \beta}{4}} e^{-\frac{2\beta}{3} x^{3/2}},
	\end{equation}
	for sufficiently large $N \geq N_0(M)$. Moreover, for any $1 \leq x \leq M(\log N)^{1/3}$, there exists $C'=C'(M)>0$ such that
	\begin{equation}\label{tw_tail_left}
		\P \Big( N^{2/3} (\lambda_N-2)<-x \Big) \leq C'  e^{-\beta x^{3}/C_0},
	\end{equation}
	for sufficiently large $N \geq N'_0(M)$, where $C_0$ is the numerical constant in (\ref{less_precise}).
\end{theorem}

The proof is based on a Green function comparison method introduced in \cite{EYY12}. First, the tail distributions can be related to a properly chosen observable in terms of the imaginary part of the normalized trace of the Green function, see Lemma \ref{lemma2} below. Then it suffices to compare this observable with its Gaussian matrix counterpart using a continuous flow in (\ref{flow}) interpolating between the matrix $H_N$ and an independent Gaussian matrix. We show in Theorem \ref{green_comparison} that the difference along the flow for a long time is much smaller than the known small deviation estimates for the Gaussian ensembles. The  key point is that we need to  follow an  event of quite small probability along the interpolation flow, \ie
the Green function comparison method is organized such that a factor with this quite small probability is 
systematically preserved in every term, see Remark~\ref{remark} for more details.

In view of (\ref{TW_right}) and (\ref{tail}), the right-tail small deviation estimate in (\ref{tw_tail}) of the largest eigenvalue of Wigner matrices is optimal, at least in  a relatively small range of $x=O((\log N)^{2/3})$. This restriction is necessary under our finite moment conditions, since  it ensures that the Tracy-Widom tails in (\ref{tw_tail}) are polynomially small. The left-tail estimate is obtained from the upper bound in (\ref{less_precise}) and hence is less precise.

We remark that our current analysis based upon the Green function comparison method can be extended \cite{PY14} to sample covariance matrices of the form $X^*X$, where $X$ is an $M\times N$ matrix with independent entries and $M/N \rightarrow \gamma \in (0,\infty)$. 
For the special Gaussian cases these are known as the Laguerre Orthogonal Ensemble (LOE) and Laguerre Unitary Ensemble (LUE) and 
tail estimates for the largest eigenvalue can be found in \cite{LR10}, \cf (\ref{less_precise}).
A sharp upper bound for the LUE with the precise constant $1/12$ in the exponent for the left tail was obtained in \cite{BDMMZ01}, 
while the correct exponential constant $4/3$ for the right tail  can be obtained similarly using simpler arguments.
With the help of the  analogous version of Theorem \ref{green_comparison} for sample covariance matrices,  corresponding small deviation estimates can easily be extended to non-Gaussian sample covariance matrices, 
assuming that matrix entries of $X$ have uniformly bounded moments. We omit the details of these straightforward 
modifications.

{\it Notation:} We will use the following standard definition of {\it stochastic domination};   its
standard arithmetic  properties can be found in Proposition 6.5 in \cite{EY17}. 
\begin{definition}\label{definition of stochastic domination}
	Let $\mathcal{X}\equiv \mathcal{X}^{(N)}$ and $\mathcal{Y}\equiv \mathcal{Y}^{(N)}$ be two sequences of nonnegative random variables. We say~$\mathcal{Y}$ stochastically dominates~$\mathcal{X}$ if, for all (small) $\epsilon>0$ and (large)~$D>0$,
	\begin{align}\label{prec}
		\P\big(\mathcal{X}^{(N)}>N^{\epsilon} \mathcal{Y}^{(N)}\big)\le N^{-D},
	\end{align}
	for sufficiently large $N\ge N_0(\epsilon,D)$, and we write $\mathcal{X} \prec \mathcal{Y}$ or $\mathcal{X}=O_\prec(\mathcal{Y})$. 
\end{definition}

Throughout the paper, we use~$c$ and~$C$ to denote strictly positive constants that are independent of~$N$. Their values may change from line to line. For $X,Y \in \R$, we write $X \ll Y$ if there exists a small $c>0$ such that $|X| \leq N^{-c} |Y|$ for large $N$. Moreover, we write $X \sim Y$ if there exist constants $c, C>0$ such that $c |Y| \leq |X| \leq C |Y|$ for large $N$. 

{\bf Acknowledgment:} We thank Zhigang Bao and Rong Ma for discussion on the existing literature.

\section{Proof of Theorem \ref{main}}

We start by introducing the notations and preliminaries that will be used in the proof.

\subsection{Local law of the Green function}
For a probability measure $\nu$ on $\R$ denote by $m_\nu$ its Stieltjes transform, \ie
\begin{align}
	m_\nu(z)\deq\int_\R\frac{\dd\nu(x)}{x-z}\,,\qquad z\in\C^+\,.
\end{align}
The Stieltjes transform of the semicircle distribution $\rho_{sc}(x):=\frac{1}{2 \pi} \sqrt{(4-x^2)_+}$ is denoted by $m_{sc}$, which is the unique solution to
\begin{equation}\label{self_eq}
	1+z m_{sc}(z)+ m_{sc}^2(z)=0\,,
\end{equation}
satisfying $\Im m_{sc}(z)>0$, for $\Im z>0$. It has the following quantitative properties, see \eg~\cite{EY17};
for any $z=E+\ii \eta$ with $|E| \leq 5$ and $0<\eta \leq  10$, we have
\begin{equation}\label{msc}
	|m_{sc}(z)| \leq 1; \qquad	{\Im}  m_{sc}(z)  \sim \begin{cases}
		\sqrt{\kappa+\eta}, & \mbox{if } E \in[-2, 2], \\
		\frac{\eta}{\sqrt{\kappa+\eta}}, & \mbox{otherwise}\,, 
	\end{cases}\quad \kappa:=\min \{ |E-2|, |E+2|  \}.
\end{equation}		

We define the resolvent or Green function of the matrix $H_N$ and its normalized trace $m_N$ by
\begin{equation}\label{Green_fun}
	G(z):=\frac{1}{ H_N-z}\,, \qquad m_N(z):=\frac{1}{N} \Tr G(z)\,, \quad\qquad z \in \C^+\,.
\end{equation}
With these notations, we now recall the following local law for the Green function $G$.
\begin{theorem}[Entry-wise local law \cite{EYY12}; Isotropic local law \cite{KY17}]\label{le theorem local law}
	For any arbitrary small $\epsilon>0$, define
	\begin{equation}\label{ddd}
		S \equiv S(\epsilon):=\big\{z=E+\ii \eta:  |E| \leq 5, N^{-1+\epsilon} \leq \eta \leq  10 \big\}.
	\end{equation}
	Then for any $z\in S$, we have
	\begin{equation}
		\max_{1 \leq i,j \leq N} | G_{ij}(z) -\delta_{ij}  m_{sc} (z) |   \prec \sqrt{ \frac{\Im  m_{sc}(z)}{N \eta}} +\frac{1}{N \eta}, \qquad | m_N(z) - m_{sc}(z) | \prec \frac{1}{N \eta}\,.
	\end{equation}
	More generally, for any deterministic unit vectors $\mathbf{v},\mathbf{w} \in \C^{N}$, we have uniformly in $z \in S$ that
	\begin{align}\label{G_0}
		\Big| \langle \mathbf{v},\big( G(z)-m_{sc}(z)I\big)  \mathbf{w} \rangle   \Big| \prec \sqrt{ \frac{\Im  m_{sc}(z)}{N \eta}} +\frac{1}{N \eta}.
	\end{align}
	
\end{theorem}

For any $1\leq j\leq N$,  we define the classical location $\gamma_j$ of the $j$-th eigenvalue $\lambda_j$ by
\begin{equation}\label{classical}
	\frac{j}{N}=\int_{-\infty}^{\gamma_j} \rho_{sc}(x) \dd x.
\end{equation}
\begin{theorem}[Rigidity of eigenvalues \cite{EYY12}]
	For any $1\leq j \leq N$, we have
	\begin{equation}\label{rigidity4}
		|\lambda_j-\gamma_j| \prec N^{-2/3} \Big( \min \{ j, N-j+1\} \Big)^{-1/3}\,.
	\end{equation}
\end{theorem}

Next, we recall  a lemma from \cite{EYY12} 
which links the tail distribution  of $\lambda_N$
to the normalized trace of the Green function $m_N$. 

\subsection{Link of the tail distribution to the Green function}

Fixing a small $\epsilon>0$, we set a truncation $E_L:=2+N^{-2/3+\epsilon}$ for the largest eigenvalue $\lambda_N$. 
Using the rigidity of eigenvalues in (\ref{rigidity4}), for any $E<E_L$, the right-tail distribution can be written as, for any large $D>0$, 
\begin{equation}
	\P \big( \lambda_N \geq E \big)=\P\big( \Tr \chi_{E}(H) \geq 1 \big)+O(N^{-D}), \qquad \chi_{E}:=\one_{[E,E_L]}.
\end{equation}
Since $\Tr \chi_{E}(H)$ takes values at integers, we have
\begin{equation}\label{right0}
	\P \big( \lambda_N \geq E \big)=\E\big[\one_{ \Tr \chi_{E}(H) \geq 1}\big]+O(N^{-D})=\E\big[F\big(\Tr \chi_{E}(H) \big) \big]+O(N^{-D}),
\end{equation}
where  $F\,:\,\R_+\longrightarrow \R_+$ is a smooth and non-decreasing cut-off function such that 
\begin{equation}\label{F_function}
	F(x)=0, \quad \mbox{if} \quad 0 \leq x \leq 1/9; \qquad F(x)=1, \quad \mbox{if} \quad x \geq 2/9,
\end{equation}
Similalry, the left-tail distribution can be written as
\begin{equation}\label{left0}
	\P \big( \lambda_N < E \big)
	=\E\big[\wt F\big(\Tr \chi_{E}(H) \big) \big]+O(N^{-D}), \qquad \wt F:=1-F.
\end{equation}

For any $\eta>0$, we define the mollifier $\theta_{\eta}$ by setting
\begin{equation}\label{le mollifier}
	\theta_{\eta}(x):=\frac{\eta}{\pi(x^2+\eta^2)}=\frac{1}{\pi} \Im \frac{1}{x-\ii \eta},
\end{equation}
and we have
\begin{equation}\label{approx}
	\Tr \chi_{E} \star \theta_{\eta}(H)=\frac{N}{\pi} \int_{E}^{E_L} \Im m_N(y+\ii \eta) \dd y\,.
\end{equation}
The following lemma in \cite{EYY12} shows that $\Tr \chi_E(H)$ can be approximated by $\Tr \chi_E\star \theta_\eta(H)$ for $\eta \ll N^{-2/3}$. 

\begin{lemma}\label{lemma2}
	Fix any small $\epsilon>0$.  Set $\eta=N^{-2/3-\epsilon}$ and $l=N^{-2/3-\epsilon/9}$.
	There exist constants $c,C>0$ such that for any $|E-2|\leq N^{-2/3+\epsilon}$,  
	\begin{equation}\label{approx3}
		\Tr \chi_{E+l} \star \theta_{\eta}(H) -N^{-\epsilon/9} \leq \Tr \chi_E(H) \leq  \Tr \chi_{E-l} \star \theta_{\eta}(H) +N^{-\epsilon/9},
	\end{equation}
	with probability bigger than $1-N^{-D}$. Thus for any large $D>0$, we have
	\begin{equation}\label{approx2}
		\E \Big[ F\Big( \Tr \chi_{E+l} \star \theta_{\eta}(H) \Big)\Big]-N^{-D} \leq \E\big[F\big(\Tr \chi_{E}(H) \big) \big]\leq   \E \Big[F\Big( \Tr \chi_{E-l} \star \theta_{\eta}(H) \Big)\Big]+N^{-D}.
	\end{equation}
\end{lemma}

We first consider the right tail distribution in Theorem \ref{main}. Fixing a large $M>0$, for any $1 \leq x \leq M(\log N)^{2/3}$, we set 
\begin{equation}\label{E}
	E_+ \equiv E_+(x)=2+N^{-2/3}x.
\end{equation}
From (\ref{right0}) and (\ref{approx2}), we have an upper bound for the right-tail distribution 
\begin{align}\label{upper_right}
	\P\Big( N^{2/3} (\lambda_N-2) \geq x \Big) \leq \E \Big[ F\Big( \Tr \chi_{E_{+}-l} \star \theta_{\eta}(H) \Big)\Big]+O(N^{-D}).
\end{align}
Combining the first inequality in (\ref{approx2}) (applied to the Gaussian ensembles and $E=E_+-2l$) and~\eqref{right0} 
with 
(\ref{tail}),  the upper bound of the right-tail distribution in (\ref{upper_right}) can be bounded by
\begin{align}
	\E^{\mathrm{G \beta E}} \Big[ F\Big( \Tr \chi_{E_{+}-l} \star \theta_{\eta}(H) \Big)\Big] \leq & \P^{\mathrm{G \beta E}}(\lambda_N \geq E_{+}-2l)+O(N^{-D})\nonumber\\
	= & \P^{\mathrm{G \beta E}} \Big( N^{2/3} (\lambda_N-2)\geq x-2N^{-\epsilon/9} \Big)+O(N^{-D})\nonumber\\
	\leq & C x^{-\frac{3 \beta}{4}}e^{-\frac{2\beta}{3} x^{3/2}},
\end{align}
where we chose $D$ sufficiently large depending on $M$. 
In general, for any fixed  integer $k \geq 1$, we have
\begin{align}\label{Gaussian}
	&\E^{\mathrm{G \beta E}} \Big[ F\Big( \Tr \chi_{E_+-kl} \star \theta_{\eta}(H) \Big)\Big]  \leq C_k x^{-\frac{3 \beta}{4}} e^{-\frac{2\beta}{3} x^{3/2}}.
\end{align}

Next we consider the left-tail distribution in Theorem \ref{main} similarly. For any $1\leq x \leq M(\log N)^{1/3}$, set 
\begin{align}
	E_- \equiv E_-(x)=2-N^{-2/3}x.
\end{align}  
From (\ref{left0}) and (\ref{approx2}), the left-tail distribution can be bounded by 
\begin{align}\label{upper_left}
	\P\Big( N^{2/3} (\lambda_N-2) < -x \Big) \leq \E \Big[ \wt F\Big( \Tr \chi_{E_{-}+l} \star \theta_{\eta}(H) \Big)\Big]+O(N^{-D}),
\end{align}
where $\wt F=1-F$. 
Using the second inequality in (\ref{approx2})  (applied the Gaussian ensembles and $E=E_-+2l$) 
with the Gaussian estimates in (\ref{less_precise}) and choosing $D$ sufficiently large depending on $M$, the upper bound for the left tail distribution in (\ref{upper_left}) can be bounded by
\begin{align}\label{Gaussian_left}
	\E^{\mathrm{G \beta E}} \Big[\wt F\Big( \Tr \chi_{E_{-}+l} \star \theta_{\eta}(H) \Big)\Big] 
	\leq & C e^{-\beta x^{3}/C_0},
\end{align}
where $C_0$ is the numerical constant in (\ref{less_precise}). 
Similar estimates for $E_-+l$ being replaced with $E_-+kl$ for any fixed interger $k\geq 1$ also hold, \cf (\ref{Gaussian}).

\subsection{Proof of Theorem \ref{main}}
To prove Theorem \ref{main} using (\ref{upper_right}) and (\ref{upper_left}), it suffices to prove the following Green function comparison theorem using the Gaussian estimates in (\ref{Gaussian}) and (\ref{Gaussian_left}).
\begin{theorem}\label{green_comparison}
	Fix a large $M>0$ and a small $\epsilon>0$. For any $1 \leq x \leq M(\log N)^{2/3}$, set $E_+=2+N^{-2/3}x$, $\eta=N^{-2/3-\epsilon}$, and $l=N^{-2/3-\epsilon/9}$. Then there exists a constant $C=C(M,\epsilon)>0$ such that
	\begin{align}\label{tail_thm}
		\Big|\E\Big[F\big(\Tr \chi_{E_+ - l} \star \theta_{\eta}(H)  \big)  \Big]-\E^{\mathrm{G \beta E}}\Big[F\big(\Tr \chi_{E_+ - l} \star \theta_{\eta}(H)  \big)  \Big] \Big| \leq C N^{-\frac{1}{6}+4\epsilon} x^{-\frac{3 \beta}{4}} e^{-\frac{2\beta}{3} x^{3/2}},
	\end{align}
	for sufficiently large $N \geq N_0(M,\epsilon)$. Moreover, for any $1\leq x\leq M(\log N)^{1/3}$ and $E_-=2-N^{-2/3}x$, there exists a constant $C'=C'(M,\epsilon)>0$ such that
	\begin{align}\label{tail_thm_2}
		\Big|\E\Big[\wt F\big(\Tr \chi_{E_- +l} \star \theta_{\eta}(H)  \big)  \Big]-\E^{\mathrm{G \beta E}}\Big[\wt F\big(\Tr \chi_{E_-+ l} \star \theta_{\eta}(H)  \big)  \Big] \Big| \leq C' N^{-\frac{1}{6}+4\epsilon} e^{-\beta x^{3}/C_0},
	\end{align}
	for sufficiently large $N \geq N'_0(M,\epsilon)$, where $C_0$ is the numerical constant in (\ref{less_precise}).
\end{theorem}

\begin{remark}\label{remark}
	We compare the above theorem to the Green function comparison in Theorem 6.3 of \cite{EYY12} 
	that yielded  an upper bound $C N^{-1/6+c\epsilon}$ for $x=O(1)$. 
	This upper bound was obtained by counting the number of off-diagonal Green function entries 
	using the local law near the edges, \ie $|G_{ij}-\delta_{ij}| \prec N^{-1/3+\epsilon}$. 
	In Theorem \ref{green_comparison} we consider the small deviation regime
	of probability of order $N^{-\frac{2\beta}{3}M}$, so just
	a few off-diagonal Green function entries are not sufficiently small to derive the desired upper bound. The key observation is that the size of the observable in terms of the function $F$ (or $\wt F$) that we compare is 
	actually much smaller as the estimates in (\ref{Gaussian}) and (\ref{Gaussian_left}) indicate. 
	We will  keep tracking the quantities related to the function $F$ (or $\wt F$)
	in every term of the Green function expansion to ensure  an upper bound 
	with an additional factor $x^{-\frac{3 \beta}{4}} e^{-\frac{2\beta}{3} x^{3/2}}$  (or $e^{-\beta x^{3}/C_0}$ for the left-tail estimate).
	An additional complication is that this smallness is explicit only for the Gaussian case,
	so we will need to transfer it to the Wigner case by a Gronwall argument. 
\end{remark}

\begin{proof}
	For simplicity, we only consider the real case $\beta=1$; 
	the complex case $\beta=2$ can be proved analogously. Moreover, we will only prove the first estimate (\ref{tail_thm}) related to the right tail distribution, the other one can be obtained similarly. To simplify the notation, we set $E \equiv E_+$.
	
	We introduce the following matrix flow interpolating between the Wigner matrix $H_N=:(h^{(0)}_{ab})$ and its Gaussian counterpart,
	$W_N$,  \ie
	\begin{equation}\label{flow}
		H(t)=\mathrm{e}^{-\frac{t}{2}}H_N +\sqrt{1-\mathrm{e}^{-t}} W_N\,,\qquad  t \in \R^+,
	\end{equation}
	where $W_N=:(w_{ab})$ is a GOE matrix defined in (\ref{GOE}) which is independent of $H_N$. We define the time-dependent Green function of the matrix $H \equiv H(t)$ by
	\begin{equation}
		G \equiv G(t,z):=\frac{1}{H(t)-z}, \qquad m_N \equiv m_N(t,z)=\frac{1}{N}\sum_{v=1}^N G_{vv}(t,z).
	\end{equation}
	In view of (\ref{approx}), we set
	\begin{equation}\label{Xidef}
		\X \equiv \X(t):=\int_{E_1}^{E_2} N \Im m_N(t,y+\ii \eta) \dd y, \qquad E_1:=E- l, \quad E_2:=2+N^{-2/3+\epsilon}.
	\end{equation}
	Recall the local laws in Theorem \ref{le theorem local law} and (\ref{msc}). Since $H(t)$ is continuous in time, we obtain the following local law for $G(t,z)$,
	\begin{equation}\label{G}
		\max_{i,j} | G_{ij}(t,z) -\delta_{ij}  m_{sc} (z) |\prec N^{-1/3+\epsilon}=:\Psi, \qquad \Im m_{sc}(z) =O(N^{-1/3-\epsilon}),
	\end{equation}
	uniformly for any $t \in \R^+$ and any $z \in S$ with $|\Re z-2| \leq N^{-2/3+\epsilon}$ and $\Im z=N^{-2/3-\epsilon}$. 
	
	Taking the time derivative of $\E[F(\X)]$ and using the following rules
	\begin{equation}\label{rule_1}
		\frac{\partial G_{ij}}{ \partial h_{ab}}=-\frac{G_{ia} G_{bj}+G_{ib} G_{aj}}{1+\delta_{ab}}, \qquad G_{ij}=G_{ji},
	\end{equation}
	we have
	\begin{align}\label{step00}
		\frac{\dd}{\dd t}\E [F(\X)]=& \E\Big[F'(\X) \sum_{v=1}^N \Im \int_{E_1}^{E_2}   \sum_{a\leq b}\dot{h}_{ab}(t)\frac{\partial G_{vv}(y+\ii \eta)}{\partial h_{ab}} \dd y  \Big]\nonumber\\
		=&-\E\Big[F'(\X) \Im \int_{E_1}^{E_2}   \sum_{a,b=1}^N \dot{h}_{ab}(t)\big[ G^2(y+\ii \eta) \big]_{ab} \dd y  \Big]\nonumber\\
		=&-\E\Big[F'(\X)   \sum_{a,b=1}^N \dot{h}_{ab}(t) 	\Delta \Im G_{ab} \Big],
	\end{align}
	where in the last step we used that $\frac{\dd G(z)}{\dd z}=G^2(z)$ and the abbreviation 
	\begin{align}\label{dim}
		\Delta \Im G \equiv (\Delta \Im G)(t,E_1,E_2):=\Im G(t,E_2+\ii \eta)-\Im G(t,E_1+\ii \eta).
	\end{align}
	
	We use the cumulant expansion formula
	(see, e.g. Lemma 3.1 in \cite{HK17}). 
	
	\begin{lemma}\label{cumulant}
		Define the $k$-th cumulant of a real-valued random variable $h$ to be
		\begin{align}\label{cumulant_k}
			c^{(k)}(h):=(-\ii)^{k} \Big(\frac{\dd^k}{\dd t^k} \log \E e^{\ii t h} \Big)\Big|_{t=0}.
		\end{align}
		Let $f: \R \longrightarrow \C$ be a smooth function and denote by $f^{(k)}$ its $k$-th derivative. Then for any fixed $l \in \N$, we have
		\begin{align}\label{cumulform}
			\E \big[h f(h)\big]=\sum_{k+1=1}^l \frac{1}{k!} c^{(k+1)}\E[ f^{(k)}(h) ]+R_{l+1}\,,
		\end{align}
		where the error term satisfies
		\begin{align}\label{error}
			|R_{l+1}| \leq C_l \E |h|^{l+1} \sup_{|x| \leq K} |f^{(l)}(x)| +C_l \E \Big[ |h|^{l+1} 1_{|h|>K}\Big] \|f^{(l)}\|_{\infty},
		\end{align}
		and $K>0$ is an arbitrary fixed cutoff.
	\end{lemma}
	
	Note that from (\ref{flow}) we have
	\begin{align}\label{relation}
		\dot h_{ab}(t)=-\frac{e^{-\frac{t}{2}}}{2} h^{(0)}_{ab}+\frac{e^{-t}}{2\sqrt{1-e^{-t}}} w_{ab},\quad  \frac{\partial G_{ij}}{ \partial h^{(0)}_{ab}}=e^{-\frac{t}{2}}\frac{\partial G_{ij}}{\partial h_{ab}}, \quad \frac{\partial G_{ij}}{ \partial w_{ab}}=\sqrt{1-\mathrm{e}^{-t}}\frac{\partial G_{ij}}{\partial h_{ab}},
	\end{align}
	where $\{h^{(0)}_{ab}\}$ and $\{w_{ab}\}$ are independent, $c^{(2)}(h^{(0)}_{ab})=c^{(2)}(w_{ab})=N^{-1}~(a\neq b)$, $c^{(2)}(w_{aa})=2N^{-1}$, $c^{(2)}(h_{aa})=m_2N^{-1}$ with $m_2$ given in (\ref{second_moment}), and $c^{(k)}(w_{ab}) \equiv 0$ for $k \geq 3$. Applying~\eqref{cumulform} on the right side of (\ref{step00})
	with $h=\dot h_{ab}(t)$ given in (\ref{relation}), since the second cumulant of the off-diagonal entries $h^{(0)}_{ab}$ and $w_{ab}$ are matched, we observe the precise cancellations of the resulting second order expansion terms and obtain that
	\begin{align}\label{step_1}
		\frac{\dd}{\dd t}\E [F(\X)] =&  \frac{1}{N} (m_2-2)\sum_{a=1}^N \E \Big[\frac{ \partial \big(F'(\X) \Dim G_{aa} \big)}{\partial h_{aa}}\Big]\nonumber\\
		&+ \sum_{k+1=3}^{2D} \frac{1}{k!} \frac{1}{N^{(k+1)/2} } \sum_{a, b=1}^N  s^{(k+1)}_{ab}(t) \E \Big[\frac{ \partial^{k} \big(F'(\X) \Dim G_{ab}\big)}{\partial h^{k}_{ab}}\Big]  +O(N^{-D})\nonumber\\
		=:&I_2+\sum_{k+1=3}^{2D} I_{k+1}+O(N^{-D}),
	\end{align}
	where
	\begin{align}\label{S_ab}
		s^{(k+1)}_{ab}(t)=e^{-\frac{(k+1)t}{2}} c^{(k+1)}(\sqrt{N} h^{(0)}_{ab}), \qquad k+1 \geq 3
	\end{align}
	with $c^{(k+1)}(\sqrt{N} h^{(0)}_{ab})$ are the $(k+1)$-th cumulants defined in (\ref{cumulant_k}) of the normalized matrix entries of the initial matrix $H_N$. Here we truncated
	the cumulant expansions at the $2D$-th order and the last error term follows from the local law in (\ref{G}) and the moment condition in (\ref{moment_condition}).

	Using the differentiation rule in (\ref{rule_1}), for any $k \in \N$, we have
	\begin{align}\label{rule_2}
		\frac{ \partial F^{(k)}(\X) }{ \partial h_{ab} }=-\frac{2F^{(k+1)}(\X)}{1+\delta_{ab}}  \sum_{v=1}^N \Im  \Big(\int_{E_1}^{E_2} G_{va} G_{bv} \dd y \Big)=-\frac{2}{1+\delta_{ab}} F^{(k+1)}(\X) \Delta \Im G_{ab}\,,
	\end{align}
	with $\Dim$ defined in (\ref{dim}). Then the second order term on the right side of (\ref{step_1}) can be written as
	\begin{align}\label{second}
		I_2=-\frac{1}{N} (m_2-2)\sum_{a=1}^N \E \Big[ F'(\X) \Dim (G_{aa})^2 +F''(\X) (\Dim G_{aa})^2 \Big].
	\end{align}
	Using the definition of $\Dim$ in (\ref{dim}) and the local law in (\ref{G}), we have
	\begin{align}\label{I_2}
		|I_2| \prec N^{-1/3+\epsilon} \E |F'(\X)|+N^{-2/3+2\epsilon}\E |F''(\X)|. 
	\end{align}
	From the definition of the function $F$ in (\ref{F_function}), all the derivatives of $F$ are uniformly bounded, \ie for any $k\in \N$, there exist $C_k$ such that
	\begin{align}\label{F_k_C_k}
		\sup_{x\in \R}|F^{(k)}(x)| \leq C_k.
	\end{align}
	Moreover, the non-vanishing contributions to $F^{(k)}(\X)$ come from the event where
	$\X= \Tr \chi_{E -l} \star \theta_{\eta}(H)\ \in [1/9,2/9]$, thus 
	\begin{align}
		\E\big[|F^{(k)}(\X)|\big] 
		\leq C_k \P \Big( \Tr \chi_{E -l} \star \theta_{\eta}(H)\ \in [1/9,2/9]\Big).
	\end{align}
	Recall that the inequalities in (\ref{approx3}) and the rigidity of eigenvalues in (\ref{rigidity4}) imply that 
	\begin{align}
		\# \{j: \lambda_j \geq E\}-N^{-\epsilon/9} \leq \Tr \chi_{E-l} \star \theta_{\eta}(H) \leq \# \{j: \lambda_j \geq E-2l\}+N^{-\epsilon/9},
	\end{align}
	with probability bigger than $1-N^{-D}$. If $\Tr \chi_{E - l} \star \theta_{\eta}(H)\ \in [1/9,2/9]$, 
	then the largest eigenvalue $\lambda_N$
	lies  in $[E-2l,E]$. 
	Thus we have
	\begin{align}\label{F_k}
		\E\big[|F^{(k)}(\X)|\big] \leq& C_k \P\big(\lambda_N  \in[ E-2l,E]\big)+O(N^{-D})\nonumber\\
		\leq& C_k\E \Big[F\Big( \Tr \chi_{E-3l} \star \theta_{\eta}(H) \Big)\Big]+O(N^{-D}),
	\end{align}
	where in the last step we used (\ref{right0}) and the second inequality in (\ref{approx2}). Therefore, combining (\ref{F_k}) with (\ref{I_2}), we have
	\begin{align}\label{I_22}
		|I_2|\prec N^{-1/3+\epsilon}  \E \Big[F\Big( \Tr \chi_{E-3l} \star \theta_{\eta}(H) \Big)\Big]+O(N^{-D}).
	\end{align}

	Similarly, using the differentiation rules in (\ref{rule_1}) and (\ref{rule_2}), the third order term on the right side of (\ref{step_1}) can be written as a linear combination of the following terms (we ignore the deviation caused by differentiating with diagonal entries in (\ref{rule_1}) and (\ref{rule_2}) for simplicity)
	\begin{align}\label{third}
		&\frac{1}{N^{\frac{3}{2}} }  \sum_{a, b} s^{(3)}_{ab}(t) \E \Big[ F'(\X) \Dim \big( G_{aa}G_{bb} G_{ab}\big)\Big],\quad \frac{1}{N^{\frac{3}{2}} }  \sum_{a, b} s^{(3)}_{ab}(t) \E \Big[ F'(\X) \Dim \big( (G_{ab})^3\big)\Big],\nonumber\\
		& 	\frac{1}{N^{\frac{3}{2}} }  \sum_{a, b} s^{(3)}_{ab}(t) \E \Big[ F''(\X) \Dim G_{ab} \Dim \big( G_{aa}G_{bb}\big)\Big],\quad 	\frac{1}{N^{\frac{3}{2}} } \sum_{a, b} s^{(3)}_{ab}(t) \E \Big[ F''(\X) \Dim G_{ab} \Dim \big( (G_{ab})^2\big)\Big],\nonumber\\
		& \frac{1}{N^{\frac{3}{2}} } \sum_{a, b} s^{(3)}_{ab}(t) \E \Big[ F'''(\X) \Dim G_{ab} \Dim G_{ab} \Dim G_{ab}\Big].
	\end{align}
	For example, the first term above can be estimated using the local law in (\ref{G}) as
	\begin{align}
		& \frac{1}{N^{3/2} }  \sum_{a, b} s^{(3)}_{ab}(t) \E \Big[ F'(\X) \Dim \big( G_{aa}G_{bb} G_{ab}\big)\Big]  \nonumber\\
		=&\frac{1}{N^{3/2} }  \sum_{a, b} s^{(3)}_{ab}(t) \E \Big[ F'(\X) \Dim \big( m^2_{sc} G_{ab}\big)\Big]+O_{\prec}(\sqrt{N}\Psi^2 \E|F'(\X)|)\nonumber\\
		=&O_{\prec}\big((\Psi+\sqrt{N}\Psi^2) \E|F'(\X)|\big),
	\end{align}
	where in the last step we used the isotropic local law in (\ref{G_0}) with $\mathbf{w}=N^{-\frac{1}{2}}\big(s^{(3)}_{a1}(t),s^{(3)}_{a2}(t),\cdots,s^{(3)}_{aN}(t)\big)$ and $\mathbf{v}=\mathbf{e}_{a}$, \ie
	\begin{align}
		\frac{1}{\sqrt{N}} \sum_{b=1}^N s^{(3)}_{ab}(t) G_{ab} \prec N^{-1/3+\epsilon}=\Psi,
	\end{align}
	uniformly for any $t \in \R^+$ and any $z \in S$ with $ |\Re z-2| \leq N^{-2/3+\epsilon}$ and $\Im z=N^{-2/3-\epsilon}$. The other terms in (\ref{third}) can be estimated similarly by power counting using the local law in (\ref{G}). Thus the third order term can be bounded by
	\begin{align}
		|	I_3| \prec& N^{-1/6+2\epsilon} \E\big[|F'(\X)|+|F''(\X)|+|F'''(\X)|\big].
	\end{align}
	Combining with the estimates in (\ref{F_k}), we have
	\begin{align}\label{I_3}
		|	I_3|	\prec& N^{-1/6+2\epsilon}  \E \Big[F\Big( \Tr \chi_{E-3l} \star \theta_{\eta}(H) \Big)\Big]+O(N^{-D}).
	\end{align}
	
	In general, for any fixed $k+1 \geq 4$, the term $I_{k+1}$ on the right side of (\ref{step_1}) is given by
	\begin{align}
		I_{k+1}= \frac{1}{N^{(k+1)/2} } \sum_{a, b=1}^N   \frac{s^{(k+1)}_{ab}(t)}{k!} \sum_{s=0}^{k} \binom{k}{s} \E \Big[\frac{ \partial^{s} \big(F'(\X)\big)}{\partial h^{s}_{ab}} \frac{ \partial^{k-s} \big( \Dim G_{ab}\big)}{\partial h^{k-s}_{ab}}\Big].
	\end{align}
	Since $\partial/\partial h_{ab}$ and the operation $\Dim$ in (\ref{dim}) commute, using the differentiation rules in (\ref{rule_1}) and (\ref{rule_2}), each resulting term can be written in the following form
	\begin{align}\label{form}
		\frac{1}{N^{(k+1)/2} } \sum_{a, b=1}^N   \frac{s^{(k+1)}_{ab}(t)}{k!} \E \Big[ F^{(l)}(\X) \prod_{i=1}^{l} \Dim \big(\prod_{j=1}^{n_i} G_{x^{(i)}_j y^{(i)}_j}\big)\Big],
	\end{align}
	where $l$ is an interger $1 \leq l \leq k+1$,  $\sum_{i=1}^{l} n_i=k+1$, and each $x^{(i)}_j$ and $y^{(i)}_j$ represents the index $a$ or $b$. Using the local law in (\ref{G}), we have
	\begin{align}\label{im_small}
		\Big|\Im \Big(\prod_{j=1}^{n_j} G_{x^{(i)}_j y^{(i)}_j} \Big) \Big| \prec \big|\Im G_{aa}+\Im G_{bb}+\Im G_{ab}\big| \prec N^{-1/3+\epsilon}.
	\end{align}
	Using the definition of $\Dim$ in (\ref{dim}) and the estimates in (\ref{F_k}), we have from (\ref{form}) that
	\begin{align}\label{I_k+1}
		|I_{k+1}| \prec& N^{-\frac{k-3}{2}-1/3+\epsilon} \E\big[ \sum_{l=1}^{k+1} |F^{(l)}(\X)|\big]\nonumber\\
		\prec& N^{-\frac{k-3}{2}-1/3+\epsilon} \E \Big[F\Big( \Tr \chi_{E-3l} \star \theta_{\eta}(H) \Big)\Big]+O(N^{-D}).
	\end{align}
	In fact, the additional  $N^{-1/3+\epsilon}$ factor coming  from considering the imaginary part
	in (\ref{im_small}) is necessary only for $k+1=4$; for  $k+1>4$,  the factor $N^{-\frac{k-3}{2}}$ 
	alone would already be sufficient.

	Therefore, combining (\ref{I_22}), (\ref{I_3}), (\ref{I_k+1}) with (\ref{step_1}), we have
	\begin{align}\label{step1}
		\frac{\dd}{\dd t}\E \Big[ F\Big(\Tr \chi_{E-l} \star \theta_{\eta}(H(t)) \Big) \Big] \leq  N^{-1/6+3\epsilon}\E \Big[ F\Big(\Tr \chi_{E-3l} \star \theta_{\eta}(H(t)) \Big) \Big]+O(N^{-D}).
	\end{align}
	Since the function $F$ in (\ref{F_function}) is uniformly bounded by one, we integrate (\ref{step1}) over $[0,T]$ where $T:= 100 M^{3/2} \log N$ with $M$ being fixed as given in Theorem \ref{green_comparison} and obtain
	\begin{align}\label{compare_1}
		\Big|\E \Big[ F\Big(\Tr \chi_{E-l} \star \theta_{\eta}(H(0)) \Big) \Big]-\E \Big[ F\Big(\Tr \chi_{E-l} \star \theta_{\eta}(H(T)) \Big) \Big]\Big|=O(TN^{-1/6+3\epsilon}).
	\end{align}
	From the definition of $H(t)$ in (\ref{flow}) and the moment condition in (\ref{moment_condition}), we have $\|H(T)-W\|_{\mathrm{max}} \prec N^{-50M^{3/2}}$ with $W=H(\infty)$ being the GOE matrix. Here $\|A\|_{\mathrm{max}}:= \max_{ij}|A_{ij}|$ for any $A \in \C^{N\times N}$. 
	Using  $\|A\|_{\mathrm{max}} \leq \|A\| \leq N \|A\|_{\mathrm{max}}$ and that $\|G(E+\ii \eta)\| \leq \frac{1}{\eta}$, $G(T,z)$ is close to the Green function of the GOE matrix $W$, denoted by $G^{W}$, \ie
	\begin{align}\label{approxxxx}
		\|G(T,z) -G^{W}(z)\|_{\mathrm{max}} \leq &\|G(T,z)(H(T)-W)G^{W}(z)\| \nonumber\\
		\leq& \frac{N}{\eta^2} \|(H(T))-W\|_{\mathrm{max}} \prec \frac{N^{7/3+2\epsilon}}{N^{50M^{3/2}}}.
	\end{align}
	Since $F$ is a smooth function with uniformly bounded derivatives, we have
	\begin{align}\label{sandwich2}
		\Big|\E \Big[ F\Big(\Tr \chi_{E-l} \star \theta_{\eta}(H(T)) \Big) \Big]-\E^{\mathrm{GOE}} \Big[ F\Big(\Tr \chi_{E-l} \star \theta_{\eta}(H) \Big) \Big]\Big| \prec \frac{N^{8/3+3\epsilon}}{N^{50M^{3/2}}}.
	\end{align}

	Therefore, combining (\ref{compare_1}) with (\ref{sandwich2}), we obtain an initial estimate
	\begin{align}\label{compare}
		&\Big|\E \Big[ F\Big(\Tr \chi_{E-l} \star \theta_{\eta}(H(0)) \Big) \Big]-\E^{\mathrm{GOE}} \Big[ F\Big(\Tr \chi_{E-l} \star \theta_{\eta}(H) \Big) \Big]\Big|\nonumber\\
		&\qquad \qquad=O\big(TN^{-1/6+3\epsilon}+\frac{N^{8/3+4\epsilon}}{N^{50M^{3/2}}}\big)=:\mathcal{E}_0,
	\end{align}
	where the second error term is much smaller than the target error in (\ref{tail_thm}) for $1\leq x\leq M (\log N)^{2/3}$. 
	
	We remark that the estimate in (\ref{step1}) holds true if we replace $E-l$ by $E-(2k+1)l$ for any fixed $k\in \N$. In general, for any fixed $k \in \N$, we have
	\begin{align}\label{stepk}
		\frac{\dd}{\dd t}\E \Big[ F\Big(\Tr \chi_{E-(2k+1)l} \star \theta_{\eta}(H) \Big) \Big] \leq  N^{-\frac{1}{6}+3\epsilon}\E \Big[ F\Big(\Tr \chi_{E-(2k+3)l} \star \theta_{\eta}(H) \Big) \Big]+O(N^{-D}).
	\end{align}
	For $k=1$, integrating (\ref{stepk}) as in (\ref{compare_1})-(\ref{compare}) and using the Gaussian estimate in (\ref{Gaussian}), we obtain that
	\begin{align}\label{step2}
		\E \Big[ F\Big(\Tr \chi_{E- 3l} \star \theta_{\eta}(H) \Big) \Big]=O\big(\mathcal{E}_0+x^{-\frac{3}{4}}e^{-\frac{2}{3} x^{3/2}}\big).
	\end{align}
	Plugging (\ref{step2}) into (\ref{step1}), we have
	$$\Big|\frac{\dd}{\dd t}\E \Big[ F\Big(\Tr \chi_{E-l} \star \theta_{\eta}(H) \Big) \Big] \Big|=
	O\big(N^{-1/6+3\epsilon}\mathcal{E}_0+N^{-1/6+3\epsilon} x^{-\frac{3}{4}} e^{-\frac{2}{3} x^{3/2}}\big),$$
	which further implies as in (\ref{compare}) that
	$$\Big|\big(\E- \E^{\mathrm{GOE}}  \big) \Big[ F\Big(\Tr \chi_{E-l} \star \theta_{\eta}(H) \Big) \Big]\Big|=O\Big(TN^{-1/6+3\epsilon}\mathcal{E}_0+T N^{-1/6+3\epsilon} x^{-\frac{3}{4}} e^{-\frac{2}{3} x^{3/2}}+\frac{N^{8/3+4\epsilon}}{N^{50M^{3/2}}}\Big).$$
	In this way, we have improved the initial estimate in (\ref{compare}) by an additional factor $TN^{-1/6+3\epsilon}$ with $T=100M^{3/2}\log N$ plus two irrelevant additive terms much smaller than the target error in (\ref{tail_thm}). We can clearly iterate this procedure to gain additional factors $TN^{-1/6+3\epsilon}$ in each step. 
	In order to obtain the target error in (\ref{tail_thm}), we start from considering (\ref{stepk}) for $k=D_0$, where $D_0$ is chosen sufficiently large so that $\big(T N^{-1/6+3\epsilon}\big)^{D_0-1} \ll  x^{-3/4} e^{-\frac{2}{3} x^{3/2}}$. Since $1\leq x \leq M (\log N)^{2/3}$ for a fixed $M>0$, $D_0$ can be
	chosen a large constant depending on $M$ and $\epsilon$. Performing the above arguments iteratively, we obtain the desired upper bound in Theorem \ref{green_comparison}.

\end{proof}

\section{Appendix: Proof of Lemma \ref{lemma_tail}}

Before giving the proof of Lemma \ref{lemma_tail}, we recall some useful properties of the Gaussian ensembles, for references see \cite{AGZ10} and \cite{M04}. In order to keep consistent with notations used before, we consider the rescaled Gaussian matrix $X_N=\sqrt{\frac{N}{2}}H$ with $H$ defined in (\ref{GOE}) and (\ref{GUE}). 
Denote the eigenvalues of the rescaled matrix $X_N$ by $(\nu_j)_{j=1}^N$ in non-decreasing order. The joint eigenvalue density is given~by
\begin{align}\label{prop}
	{p}(\nu_1, \cdots, \nu_N)=\frac{1}{Z_{N,\beta}} \prod_{i < j} |\nu_i-\nu_j|^\beta \mathrm{e}^{-\frac{\beta}{2} \sum_{i=1}^N \nu_i^2}, \qquad \beta=1,2,
\end{align}
with $Z_{N,\beta}$ being the normalization constant. For $\beta=2$, the distribution of the eigenvalues is well known to form a 
determinantal point process. Let $q_{k}$ be the $k$-th Hermite orthogonal polynomial given by
\begin{align}\label{le hermite functions}
	q_k(x):=(-1)^k \mathrm{e}^{x^2} \frac{\dd^k}{\dd x^k} \mathrm{e}^{-x^2}, \qquad 	\phi_k(x):=\frac{1}{\sqrt{2^k k! \sqrt{ \pi}}} \mathrm{e}^{-\frac{x^2}{2}} q_k(x)\,.
\end{align}
Then the $n$-point correlation function of the eigenvalue process is given by
\begin{align}
	p_n(\nu_1, \cdots, \nu_n)=\det[ K_{N,2}(\nu_i,\nu_j)]_{1 \leq i,j \leq n},
\end{align}
where the correlation kernel can be written as (for a reference see \eg \cite{A05})
\begin{align}\label{kernel}
	K_{N,2}(x,y)=&\sqrt{\frac{N}{2}} \int_0^\infty  \big[\phi_{N}(x+z) \phi_{N-1}(y+z)+\phi_{N-1}(x+z)\phi_{N}(y+z) \big]\dd z.
\end{align}
Since $X_N=\sqrt{\frac{N}{2}}H$, we have
\begin{align}\label{integral}
	\P^{\mathrm{GUE}}\big( \lambda_N > 2+N^{-2/3} r\big)=&\P\big( \nu_N > \sqrt{2N}+\frac{r}{\sqrt{2}N^{1/6}}\big) \leq \E[\#\{i:\nu_i >\sqrt{2N}+\frac{r}{\sqrt{2}N^{1/6}}\}]\nonumber\\
	=&\int^\infty_{\sqrt{2N}+\frac{r}{\sqrt{2}N^{1/6}}}  K_{N,2}(x,x) \dd x=\int_{r}^\infty  \wt K_{N,2} (\wt x,\wt x)\dd \wt x,
\end{align}
where in the last step we changed the variable $x=\sqrt{2N}+\frac{\wt x}{\sqrt{2}N^{1/6}}$ and 
the rescaled kernel is given by
\begin{align}\label{edge}
	\wt K_{N,2}( x, y):=\frac{1}{\sqrt{2} N^{1/6}} K_{N,2}\Big( \sqrt{2N}+\frac{ x}{\sqrt{2} N^{1/6}},\sqrt{2N}+\frac{ y}{\sqrt{2}N^{1/6}} \Big).
\end{align}
Using the intergral form in (\ref{kernel}), we rewrite it as
\begin{align}\label{kernel_xy}
	\wt K_{N,2}( x, y)=&\frac{1}{2\sqrt{2}} \int_0^\infty  \big[f(x+z) g(y+z)+g(x+z)f(y+z)\big] \dd z,
\end{align}
with 
\begin{align}\label{fg}
	f(x):=N^{1/12} \phi_N \big(\sqrt{2N}+\frac{x}{\sqrt{2} N^{1/6}}\big), \qquad g(x):=N^{1/12} \phi_{N-1} \big(\sqrt{2N}+\frac{x}{\sqrt{2} N^{1/6}}\big).
\end{align}
To estimate the function $f$ and $g$, we recall the uniform asymptotics of Hermite polynomials \cite{S59}
\begin{align}\label{sk}
	q_N(\sqrt{2N+1}t)=&(2 \pi )^{1/2} (2N+1)^{N/2+1/6} e^{\frac{1}{2} (2N+1)(t^2-\frac{1}{2})} \big(\frac{\xi}{t^2-1} \big)^{1/4}\mathrm{Ai}((2N+1)^{2/3} \xi)  \nonumber\\
	&\quad \times \big( 1+O(N^{-1}(1+t^2)^{-1})\big),
\end{align}
where 
\begin{align}\label{xi}
	\xi=&\big(\frac{3}{4} t \sqrt{t^2-1}-\frac{3}{4} \mathrm{arccosh} t \big)^{2/3}=2^{1/3}(t-1)+O((t-1)^2), \qquad  t\geq 1,\nonumber\\
	\xi=&\big(-\frac{3}{4} t \sqrt{1-t^2}+\frac{3}{4} \mathrm{arccos} t \big)^{2/3}=2^{1/3}(t-1)+O((t-1)^2), \qquad  t\leq 1,
\end{align}
and $\mathrm{Ai}$ is the Airy function with the following asymptotic behavior
$$\mathrm{Ai}(x) \sim \frac{e^{-\frac{2}{3} x^{3/2}}}{2 \sqrt{\pi} x^{1/4}}, \qquad x \rightarrow \infty.$$
In view of (\ref{le hermite functions}) and (\ref{fg}), we set $\sqrt{2N+1}t=\sqrt{2N}+\frac{x}{\sqrt{2} N^{1/6}}$ for any $1\leq x \leq N^{1/4}$. Then we have $t=1+\frac{x}{2N^{2/3}}+O(\frac{1}{N})$ and thus from (\ref{xi})
$$(2N+1)^{2/3} \xi=x+O\big(\frac{1}{N^{1/3}}+\frac{x^2}{N^{2/3}} \big).$$ Combining with (\ref{le hermite functions}), (\ref{fg}), and (\ref{sk}), using the Stirling approximation, we have 
\begin{align}\label{bound}
	|f(x)| \leq C \big(\frac{\xi}{t^2-1} \big)^{1/4}\mathrm{Ai}((2N+1)^{2/3} \xi) \leq C' x^{-1/4}e^{-\frac{2}{3} x^{3/2}}\Big(1+O\big(N^{-1/24}\big)\Big).
\end{align}
The function $g$ can be estimated similarly. We hence obtain from (\ref{integral}) and (\ref{kernel_xy}) that
\begin{align}\label{result_1}
	\P^{\mathrm{GUE}}\big( \lambda_N > 2+N^{-2/3} r\big)\leq &	\frac{1}{\sqrt{2}}  \int_{r}^{\infty}  \int_0^\infty  f(x+z) g(x+z) \dd z \dd x\nonumber\\
	=&\frac{1}{\sqrt{2}} \int_r^{\infty} (s-r)f(s) g(s) \dd s \leq C r^{-3/2}e^{-\frac{4}{3} r^{3/2}}.
\end{align}
This proves Lemma \ref{lemma_tail} in the $\beta=2$ case.

We next consider the GOE in a similar way. For $\beta=1$ in (\ref{prop}), the resulting eigenvalue process is a Pfaffian point process. In view of (\ref{integral}), it suffices to estimate the corresponding one-point correlation function which can be written as 
(see equation (6.3.8) in \cite{M04})
\begin{align}\label{one_point_real}
	K_{N, 1}(x,x)=&K_{N,2}(x,x)+\sqrt{\frac{N}{2}}\phi_{N-1}(x) \Big( \int_{-\infty}^{\infty} \mathrm{sgn}(x-t) \phi_N(t) \dd t \Big)+\frac{1}{2 I_{N-1}} \phi_{N-1}(x) \one_{N=2m+1},
\end{align}
where $K_{N,2}(x,x)$ is the one-point correlation function for $\beta=2$ given in (\ref{kernel}), and
\begin{align}\label{int_even}
	I_{2m}:=\int_0^{\infty} \phi_{2m}(t) \dd t = \frac{2^{\frac{1}{4}}\sqrt{(2m)!}}{\pi^{\frac{1}{4}}2^m m!}
	\sim m^{-\frac{1}{4}}, \qquad I_{2m+1}:=\int_0^{\infty} \phi_{2m+1}(t) \dd t=O(m^{-\frac{1}{4}}).
\end{align}
For notational simplicity, we only consider the case when $N$ is even, and the odd $N$ case can be estimated similarly. Since $N$ is even, the one point function can be written as
\begin{align}
	K_{N, 1}(x,x)=K_{N,2}(x,x)+\sqrt{2N} \phi_{N-1}(x) I_N -\sqrt{2N} \phi_{N-1}(x) \int_{x}^{\infty}\phi_N(t) \dd t.
\end{align}
Changing the variable as in (\ref{edge}), the rescaled one point correlation function can be written as
\begin{align}\label{K_edge_goe}
	\wt	K_{N, 1}(x,x):=\wt	K_{N, 2}(x,x)+ N^{1/4} I_{N} g(x) - \frac{1}{\sqrt{2}} g(x)\int_{x}^{\infty} f(s) \dd s,
\end{align}
with $f$ and $g$ given in (\ref{fg}). Combining with the estimates in (\ref{bound}), (\ref{result_1}) and (\ref{int_even}), we have
\begin{align}
	\P^{\mathrm{GOE}}\big( \lambda_N > 2+N^{-2/3} r\big) \leq \int_{r}^{\infty} \wt	K_{N, 1}(x,x) \dd x \leq C x^{-3/4} e^{-\frac{2}{3} r^{3/2}}.
\end{align}
We hence finished the proof of Lemma \ref{lemma_tail}.

\end{document}